\documentclass[a4paper]{article}

\usepackage[english]{babel}
\usepackage[utf8x]{inputenc}

\usepackage{amsmath}
\usepackage{graphicx}

\usepackage{epsfig} % for postscript graphics files
\usepackage{epstopdf} % converts eps fig.s to pdf fig.s for PDFTexify
\usepackage{times} % assumes new font selection scheme installed
\usepackage{amsmath} % assumes amsmath package installed
\usepackage{amssymb}  % assumes amsmath package installed
\usepackage{amsfonts}
\usepackage{cite}

\newtheorem{theorem}{Theorem}[section]

\title{Role of Iso-connectivity Topologies in Multi-agent Interactions}

\author{Rajdeep~Dutta,~and~Daniel~Pack% <-this % stops a space
\thanks{R. Dutta is with the Unmanned Systems Laboratory,  Department of Electrical Engineering,
University of Texas at San Antonio, TX-78249, USA (e-mail: rajdeep.dutta@utsa.edu).}% <-this % stops a space
\thanks{D. Pack is with the College of Engineering and Computer Science, University of Tennessee at Chattanooga, TN-37403, USA (daniel-pack@utc.edu).}}

\begin{document}
\maketitle

%\date{March 22, 2016}

\begin{abstract}
In this paper, we present the benefits of exploring different topologies with equal connectivity measure, or \textit{iso-connectivity topologies}, in relation to the multi-agent system dynamics. The level of global information sharing ability among agents in a multi-agent network can be quantified by a connectivity measure, called as \textit{the Algebraic Connectivity} of the associated graph consisting of point-mass agents as nodes and inter-connection links between them as edges. Distinct topologies with the same connectivity play profound role in multi-agent dynamics as they offer various ways to reorganize agents locations according to the requirement during a cooperative mission, without sacrificing the information exchange capability of the entire network. Determination of the distinct multi-agent graphs with identical connectivity is a multimodal problem, in other words, there exist multiple graphs that share the same connectivity. We present analytical solutions of determining distinct graphs with identical connectivity. A family of isospectral graphs are found out by utilizing an appropriate similarity transformation. Moreover, a zone of no connectivity change in a dense graph is unraveled where an agent can move freely without causing any change in the global connectivity. The proposed solutions are validated with the help of sufficient examples.
\end{abstract}

{ \small
\noindent  {\bf Mathematics Subject Classification (2010).}  ~~05C40, 05C50, 05C62, 93A30, 94C15.
\\
\noindent  {\bf  Keywords.} ~~ Multi-agent graphs, algebraic connectivity, matrix algebra, calculus of variation.
}

\maketitle


\bibliographystyle{amsplain}

\section{Introduction}
Realization of multi-agent topologies facilitates in understanding network characteristics, agents interactions and group dynamics in wide range of engineering applications as well as in pedagogical courses \cite{overview, ref_LMP}. In network control and dynamics research community, the analysis and synthesis of multi-agent topologies have drawn serious attention over the past \cite{overview, WRen, Automatica}.  Most of the real world network topologies can be represented mathematically by using their corresponding graphs. For instance, in cooperative robotics applications \cite{overview}, a multi-agent system is represented by a graph where nodes are the point-robots and edges are the inter-connections between them. In a decentralized network of multiple agents, the group performance and task accomplishment depend on agents information sharing ability to their neighbors. The connectivity measure of a graph, i.e. the \textit{algebraic connectivity}, plays vital role in group dynamics as it determines how well the agents can communicate to each other. A connectivity value above zero guarantees that there exists at least one spanning tree \cite{Book_GraphTh} or an information flow path among the members in a group. Earlier research shows a fair amount of contributions \cite{overview, Ghosh_Boyd, Pap_conn2} dedicated to maintaining, controlling and maximizing connectivity of a multi-agent network, during cooperative tasks such as target tracking and formation control \cite{RD_IROS14}.

Isospectral systems \cite{rspa_Dutta} exhibit similar dynamic characteristics since they resonate at the same natural frequencies. The family of graphs that share the same spectrum are called isospectral \cite{isospectral_graph}. Isospectral graphs also have identical connectivity, however, the converse is not true always. This is because the connectivity of a graph can be measured just by computing one eigenvalue of the corresponding Laplacian; the entire spectrum is not required. There exist different multi-agent graphs that have the same connectivity measure, and so the inverse problem of identifying a particular topology from just the connectivity is not unique. The present work focuses on finding distinct graphs with the same connectivity. These graphs \cite{isograph_generation} have significant applications in topology categorization and reconfiguration \cite{reconfig} during a cooperative mission. Topology categorization offers a variety of formation options to the agents in a group, which can be selected according to the specific application, constraints or surrounding environment during a cooperative task.
%Figure \ref{fig:overview} shows a group of aerial vehicles flying in different formations depending on the dynamic environment.
%\begin{figure}[h]
%  \centering
%   \includegraphics[width=11cm,height=6cm]{Figs/overview_formation}\\
%  \caption{A group of aerial vehicles flying over dynamic environment in different formations.}\label{fig:overview}
%\end{figure}

The relationship between a multi-agent system dynamics and its corresponding graph properties is critical in terms of the network controllability and stability \cite{siam_controllability}, \cite{siam_stabilization}. For a multi-agent network with one(multiple) leader(s), Rahmani et. al. established how the symmetry structure of a network characterized by automorphism group, relates to the controllability of such network. Tanner \cite{cnvty_Cntrllability_Tanner} investigated the effects of network connectivity and size on its dynamic controllability in presence of a single leader, considering neighborhood based inter-connection topology structure. In multi-agent cooperative tasks, an increase in the network connectivity helps in faster convergence of the consensus \cite{WRen}. The author in \cite{cnvty_Cntrllability_Tanner} derived sufficient conditions for the overall network controllability, and also showed that high connectivity may cause adverse effect on network controllability. In \cite{cnvty_Cntrllability3}, the structural controllability problem of multi-agent network with single leader was investigated under switching topologies.

The rest of the paper is organized as follows. In Section 2, we brief the mathematical representation of a multi-agent graph and state the problem objective. Section 3 presents an analytical method of generating graphs with identical spectrum using the similarity transformation. In Section 4, we describe how to find graphs with the same connectivity in a network with just one mobile agent. Finally, Section 5 concludes the present work.

\section{Problem Formulation}

\subsection{Graph Theory Preliminaries}
A graph $G$ consisting of $n$ agents can be captured mathematically \cite{Book_GraphTh} by the associated Laplacian matrix $L_n(G)$ of order $n$. In this work, all the graphs are assumed as undirected with bidirectional inter-agent communication links.  

\noindent  The Laplacian matrix is the difference between the diagonal matrix of vertex degrees and the adjacency matrix of inter-agent connections, i.e. $L_n(G)= D_n(G) - A_n(G)$. Laplacian is a symmetric positive semi-definite matrix and its smallest eigenvalue is zero. Let, the eigenvalues of $L_n(G)$ in ascending order are $\lambda_1 < \lambda_2 \leq \lambda_3 \leq ... \leq \lambda_n $ and the corresponding eigenvectors are $\mathbf{v}_1, \mathbf{v}_2,..., \mathbf{v}_n$, respectively. The spectrum of a graph $s_n(G)$ is defined by the set of all Laplacian eigenvalues, i.e. $s_n(G) = \{\lambda_1, \lambda_2,..., \lambda_n\}$.

%----------------------
The elements of a Laplacian matrix depend on the relative distances between the corresponding pair of nodes. A state dependent Laplacian matrix can be given by
\begin{equation}\label{eqn:DminusA}
L_{n}(\textbf{x}) = D_{n}(\textbf{x}) - A_{n}(\textbf{x})~,
\end{equation}
where $\textbf{x}$ denotes the vector containing agents states, $D_{n}(\textbf{x})$ is the Degree matrix, and $A_{n}(\textbf{x})$ is the Adjacency matrix.

\noindent An weighted state dependent Adjacency matrix \cite{Pap_conn2} elements are given as
\begin{eqnarray}
& a_{ij}= \left\{
\begin{array}{lll}
& f(\| \mathbf{r}_{ij} \|) ~~~\mbox{for}~~~ i \neq j ~ ; \\
& 0 ~~~\mbox{for}~~~ i=j ~ ,
\end{array}
\right.
\end{eqnarray}
where $\mathbf{r}_{ij}$ denotes the relative distance vector between $(i,j)$, and $f(\| \mathbf{r}_{ij} \|)$ is a decreasing function with relative distance magnitude. In accordance with an exponential communication model \cite{RD_IROS14}, the function $f$ can be expressed as
\begin{equation}\label{eqn:adj}
f(\| \mathbf{r}_{ij} \|) =
\left\{
\begin{array}{lll}
& e^{- \frac{\sigma}{R} \| \mathbf{r}_{ij} \| } ~~~\mbox{if}~~~ \| \mathbf{r}_{ij} \| \leq R ~ ; \\
& 0 ~~~\mbox{if}~~~ \| \mathbf{r}_{ij} \| > R ~ .
\end{array}
\right.
\end{equation}
where $\sigma$ is a large positive constant representing the rate of decay of communication quality over distance. An inter-agent connection strength decays exponentially within the communication range $R$, and the connection is lost beyond this range. Note that every agent in the network has equal communication range.

% Note that a high constant value in equation (\ref{eqn:adj}) serves the purpose: $j \in N_i$ only if $\| \mathbf{r}_{ij} \| \leq R$.

\noindent The Degree matrix $D_n(G)$ is diagonal with elements given as
\begin{equation*}
d_i=\sum \limits_{j(\neq i)=1}^{n} a_{ij}~,
\end{equation*}

\noindent Using equation (\ref{eqn:DminusA}), the elements of Laplacian $L_{n}$ are expressed as
\begin{equation}\label{eqn:lapc}
l_{ij} =
\left\{
\begin{array}{lll}
-a_{ij} ~\mbox{for}~i\neq j; \\
\sum \limits_{j(\neq i)=1}^{n+1} a_{ij} ~\mbox{for}~ i=j
\end{array}
\right.
\end{equation}
%---------------------------------

The second smallest eigenvalue $\lambda_2$ of a Laplacian $L_n$ is known as the \textit\textit{algebraic connectivity} of the related graph of $n$ agents, and the corresponding vector is known as \textit{the Fiedler vector} \cite{Fiedler_75} . This eigenvalue provides a quantitative measure of the global communication strength of a multi-agent system. A multi-agent graph with positive $\lambda_2$ has at least one spanning tree ensuring an information flow path among agents, whereas with null $\lambda_2$ there exists no spanning tree and so the information exchange gets interrupted. A higher value of $\lambda_2$ indicates more connections among agents, which is useful in cooperative missions like multi-agent tracking \cite{WRen} and formation control \cite{RD_IROS14} for maintaining, improving or controlling time-varying connectivity during the dynamics.\\

\subsection{Objective}
The goal of the present work is to explore different multi-agent topologies that have the same connectivity measure as that of a given one. In mathematical words, given a multi-agent graph $G_{base}$ with connectivity measure $\lambda_2(G_{base})=\lambda_2^{base}$, we attempt to find out a family of distinct graphs $\{ G_i \}$ such that $\lambda_2(G_i)=\lambda_2^{base}$, where the index $i$ stands for the number of solutions.
\\

\noindent  The following presents various analytical techniques of determining distinct graphs that have the same connectivity.

\section{Generation of Isospectral Graphs}
In this section, we present an elegant way of finding isospectral graphs \cite{isograph_generation} which not only have the same algebraic connectivity but also have the entire spectrum same.

\subsection{Application of Similarity Transformation}

\begin{theorem}
Distinct graphs with identical spectrum are called isospectral, and so the corresponding Laplacian matrices associated with the graphs are isospectral as well. On applying the similarity transformation with appropriate choice of the orthonormal matrix, i.e. $L_2 = Q^T L_1 Q$ where $L_1, L_2, Q \in \Re^{n \times n}$ with $Q.\mathbf{1}=\mathbf{1}$, generates an isospectral Laplacian $L_2$ from a base Laplacian $L_1$.\\
\end{theorem}

%\begin{proof}
\noindent  \textit{Proof.} It is well known that the similarity transformation \cite{Book_HornJohn} preserves eigenvalues of an original matrix. Obviously, under such transformations applied to a base Laplacian matrix $L_1$, the newly generated matrix $L_2$ possesses the same set of eigenvalues.
\begin{equation}\label{eqn:simTrans}
L_2 = Q^T L_1 Q
\end{equation}
where $Q$ is an orthonormal matrix \cite{Book_HornJohn}, i.e. $Q^T Q = Q Q^T = I$. Note that the transformed matrix $L_2$ can be any symmetric positive semi-definite matrix which may not be of Laplacian matrix structure. In order to get a distinct Laplacian matrix on applying similarity transformation to a base Laplacian matrix, we need to search for a special class of orthonormal matrices that not only preserves the eigenvalues of a Laplacian matrix but also preserves its properties and structure.

The eigenvector $\mathbf{v}_1$ corresponding to the smallest eigenvalue $\lambda_1=0$ is along the direction $\mathbf{1}$ always, due to the fact that every row sum or column sum of a Laplacian matrix is zero. Note that the fixed eigenvector $\mathbf{v}_1$ always belongs to the null space of the Laplacian $L_1$, i.e. $L_1 . \mathbf{1} = \mathbf{0}$. Similarly, the transformed Laplacian matrix $L_2$ has to meet the following condition.
\begin{eqnarray}\label{eqn:rowsumL2}
& L_2 . \mathbf{1} = \mathbf{0} \\
& \text{or}~(Q^T L_1 Q).\mathbf{1} = \mathbf{0} \\ \label{eqn:Qcond1}
& \text{or}~ L_1 Q . \mathbf{1} = \mathbf{0}
\end{eqnarray}
An obvious solution to equation (\ref{eqn:Qcond1}) is determined by exploiting the property of a Laplacian matrix \cite{Book_GraphTh}, as follows.
\begin{equation}\label{eqn:Qcond2}
Q.\mathbf{1}=\mathbf{1}
\end{equation}
The above condition (\ref{eqn:Qcond2}) leads to the fact that in order to get a distinct Laplacian with same eigenvalues from a base Laplacian by similarity transform, the class of real orthonormal matrices $Q \neq I$ have to have every row sum or column sum as $1$. $~~~~~~~~~~~~\square$
%\end{proof}

\subsection{Example}
An example of such orthogonal matrix $Q$ is any permutation matrix $J$ because it satisfies the condition (\ref{eqn:Qcond2}). The use of permutation matrix $J$ \cite{ref1_ProcAMS} to get a different Laplacian from a base Laplacian, is described below.

Consider a group of four agents, where the corresponding Laplacian matrix of order $4$ takes shape as
\begin{equation}
  L_1= \left[
          \begin{array}{cccc}
     3  &  -1  &  -1  &  -1 \\
    -1  &   2  &  -1  &   0 \\
    -1  &  -1  &   3  &  -1 \\
    -1  &   0  &  -1  &   2
          \end{array}
        \right]~.
\end{equation}
On applying similarity transformation to the base Laplacian $L_1$ by a permutation matrix $J_1$, we get
\begin{eqnarray}
&  J_1= \left[
          \begin{array}{cccc}
            0 & 0 & 0 & 1 \\
            0 & 0 & 1 & 0 \\
            0 & 1 & 0 & 0 \\
            1 & 0 & 0 & 0 \\
          \end{array}
        \right]~,
\\
&  L_2= J_1^T L_1 J_1 = \left[
          \begin{array}{cccc}
     2  &  -1  &   0  &  -1 \\
    -1  &   3  &  -1  &  -1 \\
     0  &  -1  &   2  &  -1 \\
    -1  &  -1  &  -1  &   3
              \end{array}
        \right]~.
\end{eqnarray}
On applying similarity transformation to the same Laplacian $L_1$ by another permutation matrix $J_2$, we get
\begin{eqnarray}
&  J_2= \left[
          \begin{array}{cccc}
            1 & 0 & 0 & 0 \\
            0 & 0 & 1 & 0 \\
            0 & 1 & 0 & 0 \\
            0 & 0 & 0 & 1 \\
          \end{array}
        \right]~,\\
&  L_3= J_2^T L_1 J_2 = \left[
          \begin{array}{cccc}
     3  &  -1  &  -1  &  -1 \\
    -1  &   3  &  -1  &  -1 \\
    -1  &  -1  &   2  &   0 \\
    -1  &  -1  &   0  &   2
          \end{array}
        \right]~.
\end{eqnarray}
The above two distinct Laplacians $L_2$ and $L_3$ have the same set of eigenvalues as that of the base Laplacian $L_1$. Let, the graphs corresponding to the Laplacian matrices $L_1$, $L_2$, and $L_3$ are $G_1$, $G_2$, and $G_3$, respectively. Then, these graphs $G_1$, $G_2$, and $G_3$  possess the same spectrum and so they are isospectral.

\begin{figure}[h]
  \centering
   \includegraphics[width=11cm,height=5cm]{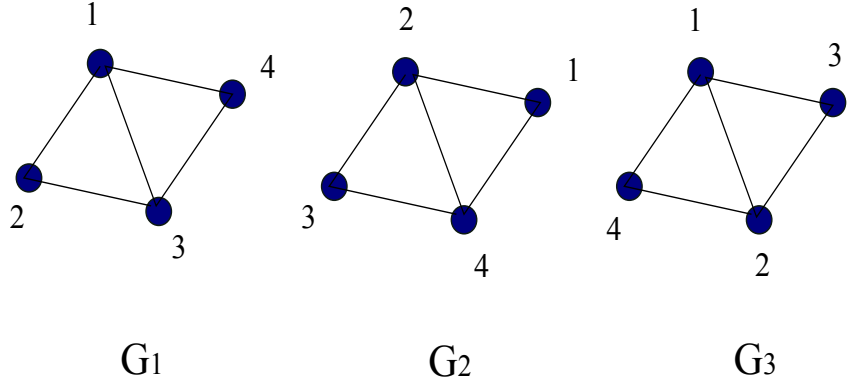}\\
  \caption{Isospectral graphs with permuted agent locations.}\label{fig:isospecGraphs}
\end{figure}

The above example shows that the rearrangement of multiple agents in a group can generate different topologies having the same spectrum. Thus, the reconfiguration of a multi-agent group by permuting agent locations has potential of switching between formations without affecting connectivity, during a cooperative task. Figure \ref{fig:isospecGraphs} shows three isospectral graphs $G_1$,$G_2$, and $G_3$, respectively.

\section{Concerns in Multi-agent Dynamics}
In the following, we explain how different multi-agent topologies with equal connectivity measure facilitate in analyzing and controlling the agent dynamics during the entire network evolves to accomplish a group mission.

\subsection{Preserving Connectivity with a Mobile Agent}

\begin{theorem}
In a dynamic network, if the differential of the corresponding Laplacian associated with the multi-agent graph becomes zero, then the network connectivity remains unchanged.
\end{theorem}

\noindent  \textit{Proof.} Consider a graph $G$ of $n$ agents which corresponds to a Laplacian matrix $L_n$. The modal equation corresponding to the second smallest eigenvalue of $L_n$ is given by
\begin{eqnarray}\label{eqn:eigenLap}
L_n \mathbf{v}_F = \lambda_2 \mathbf{v}_F~.
\end{eqnarray}
where $\lambda_2$  is the second smallest eigenvalue of $L_n$, and $\mathbf{v}_F = \mathbf{v}_2$  is the corresponding eigenvector, i.e. the Fiedler vector \cite{Fiedler_73}. The algebraic connectivity  $\lambda_2$  of $G$ can be mathematically expressed as
\begin{eqnarray}
   \lambda_2 & = & \mathbf{v}_F^T ~ L_n ~ \mathbf{v}_F ~.
\end{eqnarray}
To begin with, the graph $G$ has an algebraic connectivity of $\lambda_2^{in}$, and after moving an agent the new graph has an algebraic connectivity of $\lambda_2^{fin}$. We move the $n$th agent from its original position $P_{in}$ to new position $P_{fin}$. If all the states associated with the agent dynamics are continuous, then the change in connectivity caused by an agent movement, can be given by
\begin{eqnarray}\label{eqn:convPropg}
\lambda_2^{in} + \int_{P_{in}}^{P_{fin}} \text{d} \lambda_2 & = & \lambda_2^{fin}~.
\end{eqnarray}
%--------------------------------------
By taking differential to both sides of equation (\ref{eqn:eigenLap}), we obtain
\begin{eqnarray}\label{eqn:lam2derv}
(dL_n)\mathbf{v}_F + L_n  d\mathbf{v}_F = (d \lambda_2) \mathbf{v}_F + \lambda_2 d\mathbf{v}_F ~.
\end{eqnarray}
We now intend to determine an explicit expression of the differential of the algebraic connectivity, by pre-multiplying both sides of equation (\ref{eqn:lam2derv}) with  $\mathbf{v}_F^T$ as shown below.
\begin{eqnarray}\label{eqn:lam2derv2}
\mathbf{v}_F^T (dL_n)\mathbf{v}_F + \mathbf{v}_F^T L_n  d\mathbf{v}_F = (d \lambda_2) \mathbf{v}_F^T \mathbf{v}_F + \lambda_2 \mathbf{v}_F^T d\mathbf{v}_F
\end{eqnarray}
As the Laplacian $L_n$ is a symmetric matrix, its right and left eigenvectors do not differ \cite{Book_HornJohn}, and so the relation $\mathbf{v}_F^T L_n = \lambda_2 \mathbf{v}_F^T$ follows from equation (\ref{eqn:eigenLap}).  Also, it is to be remembered that the Fiedler vector is normalized, i.e. $\mathbf{v}_F^T \mathbf{v}_F =1$, and so equation (\ref{eqn:lam2derv2}) leads to the following expression of the differential of the algebraic connectivity.
\begin{eqnarray}\label{eqn:convDiff}
   \text{d}\lambda_2 & = & \mathbf{v}_F^T ~ \text{d}L_n ~ \mathbf{v}_F~.
\end{eqnarray}
%--------------------------------------
Thus, by using the calculus of variation \cite{eig_diff}, the differential of the algebraic connectivity can be expressed in relation to the differential of the corresponding Laplacian matrix. According to the equations (\ref{eqn:convPropg}) and (\ref{eqn:convDiff}), the connectivity $\lambda_2$ will be unchanged ($\lambda_{in}=\lambda_{fin}$) even with a mobile agent if
\begin{eqnarray}\label{eqn:diff1Lap}
 \text{d} L_n &=& 0~.~~~~~~~~~~~~~~~~~\square
\end{eqnarray}

After segregating the $n$th agent, the Laplacian $L_n$ can be expressed as follows.
\begin{eqnarray}\label{eqn:blockLap}
L_n = \left[
  \begin{array}{cc}
    L_{n-1}+B_{n-1} & -\mathbf{b}_{n-1} \\
    -\mathbf{b}_{n-1}^T & \gamma \\
  \end{array}
\right]
\end{eqnarray}
where $\gamma = \sum \limits_{j=1}^{n-1} a_{jn} \in \Re^1 $, $\mathbf{b}_{n-1}= [a_{1n}, a_{2n}, a_{3n},...]^T \in \Re^{n-1}$, $B_{n-1}=\text{diag}(\mathbf{b}_{n-1}) \in \Re^{n-1 \times n-1}$, $L_{n-1} \in \Re^{n-1 \times n-1}$, and $L_n \in \Re^{n \times n}$. It is noteworthy that the movement of $n$th agent will cause change in $\mathbf{b}_{n-1}, \gamma$ and $B_{n-1}$ usually, however, there will be no change in $L_{n-1}$. 

\subsubsection{Example}
Now, consider the $n$th agent is connected to only agents $1$ and $2$. In this case, referring to the block matrix representation (\ref{eqn:blockLap}) of Laplacian, the scalar $\gamma$, vector $\mathbf{b}_{n-1}$ and diagonal matrix $B_{n-1}$ become
\begin{eqnarray*}
& \gamma = a_{1n}+a_{2n}~, \\
& \mathbf{b}_{n-1}=[a_{1n},a_{2n},0,0,...,0]^T~, \\
& \text{and} B_{n-1} = \text{diag} (a_{1n},a_{2n},0,0,...,0)~,
\end{eqnarray*}
respectively. The Laplacian matrix $L_n$ takes shape as
\begin{eqnarray}
& L_n = \left[
  \begin{array}{cc}
    L_{n-1}+   \left(
                 \begin{array}{cccccc}
                   a_{1n} & 0 & 0 & . & . & 0 \\
                   0 & a_{2n} & 0 & . & . & 0 \\
                   0 & 0 & 0 & . & . & 0 \\
                   . & . & . & . & . & . \\
                   0 & 0 & 0 & . & . & 0\\
                 \end{array}
               \right)
     &  \left(
         \begin{array}{c}
           - a_{1n} \\
           - a_{2n} \\
           0 \\
           . \\
           0 \\
         \end{array}
       \right)
      \\ \\
~~~~~    \left(
       \begin{array}{ccccc}
         - a_{1n} & - a_{2n} & 0  & . & 0 \\
       \end{array}
     \right)
     & (a_{1n}+a_{2n}) \\
  \end{array}
\right] ~ ,
\end{eqnarray}
where $a_{1n}= \mbox{e}^{-\frac{\sigma}{R} \sqrt{(x_1 - x_n)^2 + (y_1 - y_n)^2}}$, $a_{2n}= \mbox{e}^{-\frac{\sigma}{R} \sqrt{(x_2 - x_n)^2 + (y_2 - y_n)^2}}$, and $[x_i,~y_i]~\forall ~i \in [1,n]$ denotes the $i$th agents's position in two dimensional coordinates.

\noindent  Suppose the $n$th agent moves from its original position $[x_n,~y_n]$ to a new position $[x_n^{'},~y_n^{'}]$, and the Laplacian $L_n$ corresponding to the multi-agent graph changes to $L_n^{'}$. The new Laplacian $L_n^{'}$ is given by
\begin{eqnarray}
& L_n^{'} = \left[
  \begin{array}{cc}
    L_{n-1}+   \left(
                 \begin{array}{cccccc}
                   a_{1n}^{'} & 0 & 0 & . & . & 0 \\
                   0 & a_{2n}^{'} & 0 & . & . & 0 \\
                   0 & 0 & 0 & . & . & 0 \\
                   . & . & . & . & . & . \\
                   0 & 0 & 0 & . & . & 0\\
                 \end{array}
               \right)
     &  \left(
         \begin{array}{c}
           - a_{1n}^{'} \\
           - a_{2n}^{'} \\
           0 \\
           . \\
           0 \\
         \end{array}
       \right)
      \\ \\
~~~~~    \left(
       \begin{array}{ccccc}
         - a_{1n}^{'} & - a_{2n}^{'} & 0  & . & 0 \\
       \end{array}
     \right)
     & (a_{1n}^{'}+a_{2n}^{'}) \\
  \end{array}
\right] ~ ,
\end{eqnarray}
where $a_{1n}^{'}= \mbox{e}^{-\frac{\sigma}{R} \sqrt{(x_1 - x_n^{'})^2 + (y_1 - y_n^{'})^2}}$, $a_{2n}^{'}= \mbox{e}^{-\frac{\sigma}{R} \sqrt{(x_2 - x_n^{'})^2 + (y_2 - y_n^{'})^2}}$. In order to satisfy condition (\ref{eqn:diff1Lap}), the new Laplacian elements have to be $a_{1n}^{'}=a_{1n}$ and $a_{2n}^{'}=a_{2n}$, respectively, which leads to
\begin{eqnarray}
& \| \mathbf{r}_{1n} \| = \| \mathbf{r}_{2n} \|~, \\
&~\text{or}~~  (x_1 - x_n^{'})^2 + (y_1 - y_n^{'})^2 =  (x_2 - x_n^{'})^2 + (y_2 - y_n^{'})^2 ~.
\end{eqnarray}
The above indicates that if the $n$th agent is located at any of the two intersection points of two circles, one with center $[x_1,~y_1]$ of radius $\| \mathbf{r}_{1n} \|$ and another with center $[x_2,~y_2]$ of radius $\| \mathbf{r}_{2n} \|$, then there is no change in the Laplacian $L_n$ and so the connectivity is preserved.

\begin{figure}[h]
  \centering
   \includegraphics[width=5.5cm,height=6cm]{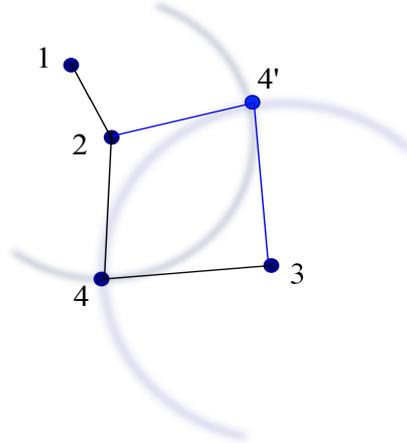}\\
  \caption{An alternative node placement without changing connectivity ; $\overline{24}=\overline{24'}$ and $\overline{34}=\overline{34'}$.}\label{fig:spclSoln}
\end{figure}

Figure \ref{fig:spclSoln} shows two graphs with identical connectivity, $G(\{ 1,2,3,4 \})$ and $G(\{ 1,2,3,4^{'} \})$, where in each graph the $4$th agent is located at either of the intersection points of the range circles of agents $2$ and $3$.

\subsection{Freedom of Moving an Agent in Dense Graphs}

\begin{theorem}
In case different Laplacian matrices associated with distinct graphs share the same algebraic connectivity and the same Fiedler vector $\mathbf{v}_F$, then that very Fiedler vector lies in the null space of the difference between corresponding Laplacians, i.e. $\Delta L . \mathbf{v}_F = \mathbf{0}$.\\
\end{theorem}

%\begin{proof}
\noindent  \textit{Proof.}  Consider two distinct graphs $G_1$ and $G_2$, which are not isospectral but have equal connectivity. Let, the Laplacian matrices corresponding to the graphs $G_1$, $G_2$ are $L_1$, $L_2$, respectively. Here, we study a specific case where different Laplacians not only share the same $2$nd smallest eigenvalue or algebraic connectivity but also share the corresponding eigenvector or Fiedler vector \cite{Fiedler_75}. The eigenvector equations of two Laplacian matrices are given by
\begin{eqnarray}\label{eqn:L1ev}
& (L_1 - \lambda_2 I ). \mathbf{v}_F = \mathbf{0}\\ \label{eqn:L2ev}
& (L_2 - \lambda_2 I ). \mathbf{v}_F = \mathbf{0}
\end{eqnarray}
Subtracting (\ref{eqn:L1ev}) from (\ref{eqn:L2ev}), we get
\begin{eqnarray}
& (L_2 - L_1). \mathbf{v}_F = \mathbf{0}~, \\ \label{eqn:diff2Lap}
& \text{or}~~\Delta L . \mathbf{v}_F = \mathbf{0}.
\end{eqnarray}
Equation (\ref{eqn:diff2Lap}) reveals that the Fiedler vector lies in the null space of the difference between two Laplacian matrices associated with two different graphs. $~~~~~~~~~~~~~~~~~~~~~~~~~~~~~\square$
%\end{proof}

Suppose the Laplacians $L_1$ and $L_2$ correspond to distinct graphs $G(t)$ and $G(t+\Delta t)$ representing a dynamic multi-agent system at time instants $t$ and $(t+ \Delta t)$, respectively. The above theorem shows the possibility of arising an unusual scenario in a group dynamics, where agent movements do not affect the network connectivity. To this end, in a dense graph with sufficient inter-agent connections, the above phenomenon may exhibit \cite{cnvty_Cntrllability_Tanner}, where the movement of one agent do not cause any change in the global connectivity. Preserving both the connectivity and the Fiedler vector in distinct graphs can be beneficial in multi-agent dynamics, which is described below using an example.

\subsubsection{Geometric Realization}
Here, we consider a graph of $4$ agents where everybody is connected to each other, and the connection strengths vary depending on the relative distances between (according to the locations of) corresponding pairs of agents. The Laplacian matrix $L_4$ associated with the graph, is given by
\begin{eqnarray}\label{eqn:L4}
L_4 = \left[
        \begin{array}{cccc}
          2+\alpha & -1 & -1 & -\alpha \\
          -1 & 2+\beta & -1 & -\beta \\
          -1 & -1 & 3 & -1 \\
          -\alpha & -\beta & -1 & 1+\alpha+\beta \\
        \end{array}
      \right].
\end{eqnarray}
where $\alpha,~\beta > 0$ are the parameters that vary during the motion of agent $4$. Note that only the last agent in the network is mobile while others are fixed.

The eigenvalues and eigenvectors of the above matrix (\ref{eqn:L4}) can be determined analytically \cite{Book_HornJohn} by solving $(L_4 - \lambda_i I) \mathbf{v}_i = 0 ~~\forall ~i \in [1,n]$, where $\lambda_i$ is the $i$th eigenvalue and $\mathbf{v}_i$ is the corresponding eigenvector. The eigenvalues of Laplacian $L_4$ are the roots of its characteristic polynomial, which is given as follows.
\begin{eqnarray}\label{eqn:polyL4}
P_{L_4}(\lambda) = \lambda (\lambda - 4) (\lambda^2 - 4 \lambda - 2 \alpha \lambda - 2 \beta \lambda + 3 + 5 \alpha + 5 \beta + 3 \alpha \beta )=0
\end{eqnarray}
By solving the characteristic polynomial (\ref{eqn:polyL4}), the Laplacian spectrum is evaluated as: $s_4 = \{0,~ 4,~ 2+\alpha+\beta \pm \sqrt{1-\alpha+\alpha^2-\beta-\alpha\beta+\beta^2} \}$. The eigenvectors of $L_4$ corresponding to the eigenvalues $0,~4$ are fixed, which are $[1,1,1,1]^T$ and $[1,1,-3,1]^T$, respectively; the other two eigenvectors are parameter dependent and change accordingly.

\noindent In order to make sure that $\lambda_2=4$ is the second smallest eigenvalue always, i.e. $\lambda_1 < \lambda_2 \leq \lambda_3 \leq \lambda_4$, the parameters $\alpha,~\beta$ have to satisfy the following inequality
\begin{eqnarray}\label{eqn:paramCond}
& 2+\alpha+\beta + \sqrt{1-\alpha+\alpha^2-\beta-\alpha\beta+\beta^2}~ > 4 .
% & \text{or}~~(1+ \alpha \beta) > (\alpha + \beta)
\end{eqnarray}
Now, we choose a set of parameters $\alpha=2, \beta=3$ and $\alpha=3, \beta=4$ satisfying (\ref{eqn:paramCond}), and then validate the claim stated above.

\noindent  The Laplacian matrix (\ref{eqn:L4}) with parameters $\alpha=2, \beta=3$ becomes
\begin{eqnarray*}
L_4^{'} = \left[
        \begin{array}{cccc}
          4 & -1 & -1 & -2 \\
          -1 & 5 & -1 & -3 \\
          -1 & -1 & 3 & -1 \\
          -2 & -3 & -1 & 6 \\
        \end{array}
      \right]=10 \times \left[
                \begin{array}{cccc}
          			 0.4 & -0.1 & -0.1 & -0.2 \\
          			-0.1 & 0.5 & -0.1 & -0.3 \\
          			-0.1 & -0.1 & 0.3 & -0.1 \\
          			-0.2 & -0.3 & -0.1 & 0.6 \\
        		\end{array}
      		  \right].
\end{eqnarray*}
The eigenvalues of the Laplacian $L_4^{'}$ are $0,~4,~5.2679,~8.7321$, and the modal matrix $M_4^{'}$ consisting of the eigenvectors of $L_4^{'}$ is given as
\begin{eqnarray*}
M_4^{'} = \left[
        \begin{array}{cccc}
          0.5000  & -0.2887  &  0.7887 &  -0.2113 \\
          0.5000  & -0.2887 &  -0.5774 &  -0.5774 \\
          0.5000  &  0.8660  &  0.0000  &  0.0000 \\
          0.5000 &  -0.2887 &  -0.2113  &  0.7887 \\
        \end{array}
      \right].
\end{eqnarray*}

\begin{figure}[h]
  \centering
   \includegraphics[width=4cm,height=5cm]{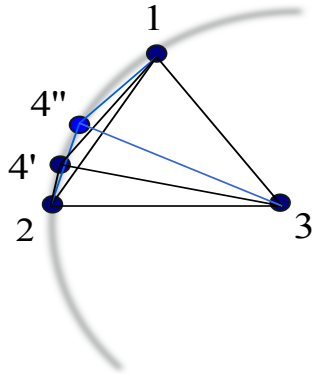}\\
  \caption{No change in connectivity zone in a graph ; $\overline{31}=\overline{32}=\overline{34'}=\overline{34''}$.}\label{fig:dense}
\end{figure}

\noindent The Laplacian matrix (\ref{eqn:L4}) with parameters $\alpha=3, \beta=4$ becomes
\begin{eqnarray*}
L_4^{''} = \left[
        \begin{array}{cccc}
          5 & -1 & -1 & -3 \\
          -1 & 6 & -1 & -4 \\
          -1 & -1 & 3 & -1 \\
          -3 & -4 & -1 & 8 \\
        \end{array}
      \right] = 10 \times \left[
        \begin{array}{cccc}
          0.5 & -0.1 & -0.1 & -0.3 \\
          -0.1 & 0.6 & -0.1 & -0.4 \\
          -0.1 & -0.1 & 0.3 & -0.1 \\
          -0.3 & -0.4 & -0.1 & 0.8 \\
        \end{array}
      \right].
\end{eqnarray*}
The eigenvalues of the Laplacian $L_4^{''}$ are $0,~4,~6.3542,~11.6458$, and the modal matrix $M_4^{''}$ consisting of the eigenvectors of $L_4^{''}$ is given as
\begin{eqnarray*}
M_4^{''} = \left[
        \begin{array}{cccc}
          0.5000  & -0.2887  &  0.7651  &  -0.2852 \\
          0.5000  & -0.2887  &  -0.6295  &  -0.5199 \\
          0.5000  &  0.8660  &  0.0000  &  0.0000 \\
          0.5000  &  -0.2887  &  -0.1355  &  0.8052 \\
        \end{array}
      \right].
\end{eqnarray*}
It can be noticed easily that the Laplacians $L_4^{'}$ and $L_4^{''}$ have the common algebraic connectivity as $\lambda_2=4$ and the common Fiedler vector as $\mathbf{v}_2 = \mathbf{v}_F = [-1, -1, 3, -1]^T$, and these matrices satisfy the condition (\ref{eqn:diff2Lap}). The corresponding graphs with same connectivity, $G(\{ 1,2,3,4^{'} \})$ and $(\{ 1,2,3,4^{''} \})$ , are shown in Figure \ref{fig:dense}. In this case, the movement of an agent doesn't influence the global connectivity which is counter intuitive \cite{siam_controllability}. In fact, these topologies offer freedom of moving an agent in a network without compromising its connectivity.

%X\newpage

\section{Conclusion}
In this work, we developed analytical methods of determining distinct multi-agent topologies with the same spectrum as well as distinct topologies with the same connectivity. We find out a family of graphs isospectral to a given one, by applying an appropriate similarity transform to the base Laplacian matrix corresponding to a given graph. A proper choice of the orthonormal matrix in the similarity transformation, leads to generate useful alternatives from a base graph while preserving connectivity. We also presented elegant strategies for preserving connectivity in a network with single mobile agent while other members are stationary. Using a simple example, we established the fact how a single agent can move freely in dense graphs without causing any change in the connectivity.

\end{document}